\def \version {2019--11--18}

\documentclass[12pt]{article}
\usepackage{amssymb}
\usepackage{amsmath}
\usepackage{fullpage}

\newcommand\F{{\cal F}}
\newcommand\K{{\cal K}}
\newcommand\cL{{\cal L}}
\newcommand\h{{\cal H}}
\newcommand\R{{\cal R}}

\newcommand\A{{\cal A}}
\newcommand\T{{\cal T}}
\newcommand\qed{\hbox{}\hfill $\Box$}

\newtheorem{thm}{Theorem}[section]
\newtheorem{prop}[thm]{Proposition}
\newtheorem{conj}[thm]{Conjecture}
\newtheorem{lem}[thm]{Lemma}

\newtheorem{cor}[thm]{Corollary}

\title{The domination number of the graph defined by two levels of the $n$-cube, II}

\author{J\'ozsef Balogh{$^{1,2,}$}\thanks{Partially supported by NSF
Grant DMS-1764123, Arnold O. Beckman Research Award (UIUC Campus Research Board RB 18132) and
the Langan Scholar Fund (UIUC).} , Gyula O.H. Katona{$^{3,}$}\thanks{The research of this author was supported
by the National Research, Development and
Innovation Office -- NKFIH Fund No's SSN117879, NK104183  and K116769.} , William Linz{$^{2}$} and Zsolt Tuza{$^{3,4,}$}\thanks{Research supported in part by the National Research, Development and Innovation Office -- NKFIH under the grant SNN 129364, and by the Sz\'echenyi 2020 grant EFOP-3.6.1-16-2016-00015.} \\ ~~ \\
\normalsize $^{1}$ Moscow Institute of Physics and Technology, Russian Federation\\
\normalsize $^{2}$ Department of Mathematics, University of Illinois at Urbana-Champaign\\
\normalsize $^{3}$ MTA R\'enyi Institute, Budapest, Hungary\\
\normalsize $^{4}$ Department of Computer Science and Systems Technology \\
\normalsize  University of Pannonia, Veszpr\'em, Hungary}

\date{\scriptsize Latest update on \version}
\begin{document}

\maketitle

\begin{abstract}
Consider all $k$-element subsets and $\ell$-element subsets $(k>\ell )$ of an $n$-element set as vertices of a bipartite graph.
Two vertices are adjacent if the corresponding $\ell$-element set is a subset of the corresponding $k$-element set.
Let  $G_{k,\ell}$ denote this graph.
The domination number of $G_{k,1}$ was exactly determined by Badakhshian, Katona and Tuza. A conjecture was also stated there on the asymptotic value ($n$ tending to infinity) of the domination number of $G_{k,2}$. Here we prove the conjecture, determining the asymptotic value of the domination number $\gamma (G_{k,2})={k+3\over 2(k-1)(k+1)}n^2+o(n^2)$.
\end{abstract}

\section{Introduction}
Let $[n]=\{ 1,2, \ldots , n\}$ be the underlying set and ${[n]\choose k}$ be the family of all $k$-element subsets of $[n]$. Suppose $n>k>\ell\geq 1$.  Consider the bipartite graph $G=\left( {[n]\choose k}, {[n]\choose \ell }; E\right)$ where the vertices $A\in {[n]\choose k}$ and
$B\in {[n]\choose \ell }$ are adjacent
 if and only if $A\supset B$. This graph will be denoted by $G_{k,\ell}$. The family ${[n]\choose k}$ is often called the $k$-th level of the $n$-cube. Then it is not much misleading to call $G_{k,\ell}$ as the graph defined by the $k$-th and $\ell$-th level.

We say that a vertex $v$ {\it dominates} the vertex $u$ in a graph $G(V,E)$ if either $u=v$ or  $\{ u,v\} \in E$. A subset $D$ of $V$ is a {\it dominating set} of the graph if every vertex $u\in V$ is  dominated by at least one element $v$ of $D$.
The {\it domination number} $\gamma (G)$ of a graph $G$ is the smallest possible size of a dominating set. Our aim is to study $\gamma(G_{k, \ell})$ for fixed $k$ and $\ell$.

Let us remark that this problem can be seen as a two-sided analogue of an old question of Erd\H{o}s and Hanani~\cite{EH}.
They defined the \textit{covering number} $M(n,k,\ell)$ to be the minimum size of a family $\K \subset \binom{[n]}{k}$ such that every $\ell$-set in $[n]$ is contained in at least one $A \in \K$. Here, in addition to covering every $\ell$-set, we also require every $k$-set to be covered. It is easy to see that dominating sets of $G_{k, \ell}$ correspond to these two-sided coverings, and that $\gamma(G_{k, \ell})$ is the two-sided ``covering number."

The value of $\gamma (G_{k,1})$ was determined by Badakhshian, Katona and Tuza in \cite{BKT}.

\begin{thm}\label{gammagk1thm} {\rm \cite{BKT}} For every\/ $k\geq 2$, we have $\gamma (G_{k,1})=n-k+1$.
\end{thm}

It seems to be much more difficult to determine $\gamma (G_{k,2})$.
A conjecture of asymptotical nature was posed in \cite{BKT}, which we prove in the present paper.

\begin{thm}\label{gammagk2thm} For every fixed\/ $k\ge 3$,\/ we have $\gamma (G_{k,2})={k+3\over 2(k-1)(k+1)}n^2+o(n^2)$
as\/ $n$ tends to infinity.
\end{thm}
 The upper bound was proved for $k=3$ and $4$ in \cite{BKT} with a construction. For general $k$, it was only proved under the assumption that certain ``small constructions" exist. In Section 2, we will prove it using a theorem of Frankl and R\"odl \cite{FR}. Of course, this proof is not constructive.

 We give two proofs for the lower bound. Both proofs are based on the {\it Graph Removal Lemma} \cite{EFR}. The first proof is much longer, but it is probably worth publishing since it gives more insight into the nature of the problem,
 showing what tools and ideas lead to different levels of approximations of the proper asymptotic value. The two proofs are in Sections 3 and 4, respectively. Section 5 contains some remarks on the case $\ell >2$, in particular a potentially tight lower bound.

\section{The upper bound}
Let $h$ be a positive integer and $\h\subset {[N]\choose h}$. The {\it degree} $d(x)$ of an element $x\in [N]$ is the number of members of $\h$ containing $x$, while the {\it two-degree} $d(\{ x,y\})$ of the pair $\{ x,y\}$ $(x,y\in [N], x\not=y)$ is the number of members of $\h$ containing both $x$ and $y$.

The proof of the upper bound is based on the following theorem.

\begin{thm}\label{frthm} {\rm \cite{FR}}
Suppose\/ $m$ is a positive integer,\/ $\epsilon> 0$ and\/ $a> 3$ are real numbers, and\/  $\h\subset {[N]\choose m}$ satisfies the following conditions. There exists a positive real\/ $\delta = \delta (\epsilon)$ such that if for some function\/
$b(N)$ one has
$$(1- \delta)b(N) < d(x) < (1 + \delta)b(N) {\mbox{\ for all\ }} x \in [N]\eqno(1)$$
and
$$d(\{x, y\}) < b(N)/(\log N)^a {\mbox{\ holds for all distinct}}\ x, y \in [N]\eqno(2)$$
 then, for all\/ $N > N_0( \delta)$, the set\/ $[N]$ can be covered by\/ ${N\over m}(1+\epsilon)$ members of\/
 $\h$.
 \end{thm}
Our task is to give a set $D$ of vertices of $G_{k,2}$ forming a dominating set and satisfying that $|D|$ is asymptotically equal to the formula given in Theorem~\ref{gammagk2thm}. The vertices of $G_{k,2}$ are two- and $k$-element subsets of $[n]$. The set of two-element subsets will be denoted by $E$, while the family of $k$-element subsets will be $\K$.
The pair $(E,\K)$ is a dominating set if and only if the following two conditions hold.
$$\mbox{ If } L\in {[n]\choose k}, L\not\in \K \mbox{ then } L \mbox{ includes an element of }E . \eqno(i)$$
 $$\mbox{ If } \{i,j\}\not\in E \mbox{ then there is a } K\in \K \mbox{ such that } \{i,j\} \subset K. \eqno(ii)$$

Suppose that $n$ is divisible by $k-1$ and form a partition $[n]=A_1\cup A_2\cup \cdots \cup A_{k-1}$ where
$|A_i|={n\over k-1}$ for every $i$. Let $E$ be the set of all pairs $\{ x,y\}$, where $x,y\in A_i$ for some $i$.
In other words, $E$ is a vertex-disjoint union of $k-1$ complete graphs on ${n\over k-1}$ vertices each. By the pigeonhole principle, every $k$-element set contains an element of $E$, therefore (i) automatically holds. Introduce the notation
$F={[n]\choose 2}-E$. This is the set of pairs connecting two distinct $A_i$s.

The family $\K$ will be constructed by applying Theorem~\ref{frthm} replacing $[N]$ by $F$. Define first the family
$$\A=\{ A:\ |A|=k, |A\cap A_i|\geq 1 {\mbox{ for every }} i\in [n]\} , $$
that is the family of all $k$-element subsets
of $[n]$ meeting every $A_i$ in exactly one element, with one exception, where the intersection has exactly two elements. A member $A\in \A$ defines a subset $H(A)$ of $F$ in the following way:
 $H(A)$ is the set of all pairs created from distinct elements of $A$, except for the only pair with both elements in the same $A_i$. Formally,
$$H(A)= \{ \{ x,y\}:\ x\not= y, x\in A, y\in A, |\{ x,y\}\cap A_i|\leq 1 {\mbox{ holds for all }} i\}.$$
The elements of $H(A)$ are the edges of a complete graph on $k$ vertices minus one edge. This is why the parameter $m$ in Theorem~\ref{frthm} will be chosen to be ${k\choose 2}-1$. Finally, define
$$\h =\{ H(A):\ A\in \A \}.$$

Theorem~\ref{frthm} will be applied for $\h$. It is easy to see that $|F|=N={k-1\choose 2}\left({n\over k-1}\right) ^2$, $m={k\choose 2}-1$.

Let us first only heuristically check the conditions (i) and (ii). The degree $d(u)$ of an element $u\in F$ should be bounded in (1). Here $u$ is a pair of elements from two distinct $A_i$'s. A member $H\in \h$ containing $u$ is determined by the set $A\subset [n]$ which contains both ends of $u$. It is obvious by symmetry that $d(u)$ does not depend on $u$ therefore this function
of $N$ satisfies (i) with say $\delta = \epsilon$. Since we have to add $k-2$ elements of $[n]$ to $u$ to obtain $A$, the order of magnitude of $d(u)$ is $n^{k-2}$. The two-degree of the pair $\{ u,v\}$ depends on whether
$u\cap v$ is empty or not. It is easy to see that its order of magnitude is smaller in the first case, in the second case it is $n^{k-3}$. Hence (ii) also holds.

To be more precise let us give the exact values of the degrees above.
$$d(u)=2\left({n\over k-1}-1\right) \left( {n\over k-1} \right)^{k-3}+(k-3){{n\over k-1}\choose 2} \left( {n\over k-1} \right)^{k-4},$$
$$d(\{ u,v\} )=\left( {n\over k-1} \right)^{k-3} {\mbox{ if $u\cap v\not= \emptyset$ and $u\cup v$ is within two $A_i$'s}},$$
$$d(\{ u,v\} )=3\left({n\over k-1}-1\right) \left( {n\over k-1} \right)^{k-4}+(k-4){{n\over k-1}\choose 2} \left( {n\over k-1} \right)^{k-5},$$
when  $u\cap v\not= \emptyset$ and $u\cup v$ meets three $A_i$'s,
$$d(\{ u,v\} )=\left( {n\over k-1} \right)^{k-4} {\mbox{ if $u\cap v= \emptyset$ and $u\cup v$ is within exactly three $A_i$'s,}}$$
and finally
$$d(\{ u,v\} )=4\left({n\over k-1}-1\right) \left( {n\over k-1} \right)^{k-5}+(k-5){{n\over k-1}\choose 2} \left( {n\over k-1} \right)^{k-6},$$
when  $u\cap v= \emptyset$ and $u\cup v$ meets four $A_i$'s.
One can see that these are in accordance with the heuristic reasoning above. The conditions of the theorem are satisfied. As a consequence, there is a subfamily $\K \subset \h$ such that its members cover every pair taken from distinct $A_i$'s and has size at most
$$|\K|={N\over m}(1+\epsilon )={{k-1\choose 2}\left({n\over k-1}\right) ^2\over {k\choose 2}-1}(1+\epsilon )\eqno(3)$$
if $N$ is large enough.

 Since every element of $F$ is covered by a member of $\K$, (ii) is also satisfied, so the pair $(E,\K)$ is really
 a dominating set in the graph $G_{k,2}$. Let us calculate its size.
 $$|E|=(k-1){{n\over k-1}\choose 2}={1\over 2(k-1)}n^2+o(n^2)\eqno(4)$$
 is obvious, and (3) implies
 $$|\K|={{k-2\over k-1}\over (k+1)(k-2)}n^2+o(n^2)={1\over (k-1)(k+1)}n^2+o(n^2).\eqno(5)$$
Adding (4) and (5) we obtain the desired asymptotic upper bound for $n$ divisible by $k-1$. It is easy to see that
the asymptotical value does not change extending the expression for other $n$. \qed

\vskip 2mm

We have given here a non-constructive, asymptotically sharp upper bound. On the other hand, in \cite{BKT} we suggested a construction with the same asymptotic behaviour which works for $k=3,4$. We encourage the reader
to check that method and try to find some construction for $k\ge 5$.

\section{Lower bound: first proof}

This proof of the lower bound is a refinement of the proof of Theorem 3 in \cite{BKT}.

It is quite obvious that our problem is closely related to the  Tur\'an theorem. For sake of completeness we repeat here
some part of \cite{BKT}.
The proof will actually be based on a stronger version of the Tur\'an theorem. Let us start with formulating the original theorem of Tur\'an. Let $T(n,s)$ denote the following graph. Partition the set $[n]$ into
 $s$ almost equal (differences of the sizes are at most one) parts: $V=V_1\cup V_2\cup \cdots \cup V_s$. Two vertices are adjacent if and only if they are in distinct $V_i$'s. The number of edges of $T(n,s)$ is denoted by $t(n,s)$.

 \begin{thm}\label{turanthm} {\rm \cite{T}} If the graph with\/ $n$ vertices contains no complete graph\/ $K_k$ as a subgraph then the number of edges cannot exceed\/ $t(n,k-1)$ and one can have equality only for\/ $T(n,k-1)$.
 \end{thm}

If one more edge is added to $T(n,k-1)$ then it creates asymptotically $\left( {n\over k-1}\right)^{k-2}$ copies of $K_k$. It is natural to guess that if
$m$ new edges are added then $m \left({n\over k-1}\right)^{k-2}$ copies of $K_k$ are obtained, if $m$ is not too large. The quantity $t(n,k-1)$
is asymptotically equal to ${n\choose 2}\left(1-{1\over k-1}\right)$. Let $m$ also be given in an asymptotic form
$m=cn$. Our above guess can be formulated in the following statement that can be easily obtained from Theorem
4 (see also Theorem 1) of a paper of Lov\'asz and Simonovits.

 \begin{cor}\label{lovaszsimcor} {\rm (of a theorem of Lov\'asz and Simonovits, \cite{LS}).} If the graph\/ $G=(V,E)$ has\/ $n$ vertices and
 $$|E|\geq {n\choose 2}\left(1-{1\over k-1}\right)+cn,$$
 where\/ $c<{1\over k-1}$, then\/ $G$ contains at least
 $$c\, {n^{k-1}\over (k-1)^{k-2}}+o(n^{k-1})$$
 copies of\/ $K_k$.
 \end{cor}

 This corollary will actually be used for the complementary graph, as follows.

 \begin{cor}\label{lovaszsimcor2}  If the graph\/ $G=(V,E)$ has\/ $n$ vertices and
 $$|E|\leq {n\choose 2}{1\over k-1}-cn,$$
 where\/ $c<{1\over k-1}$, then\/ $G$ contains at least
 $$c\, {n^{k-1}\over (k-1)^{k-2}}+o(n^{k-1})$$
 independent sets of size\/ $k$.
 \end{cor}

There the Ruzsa-Szemer\'edi theorem \cite{RSz} was used, while here the more general theorem of Erd\H os, Frankl and
R\"odl will be needed, which is a consequence of the Regularity Lemma of Szemer\'edi \cite{Sz}.

\begin{thm}\label{remlemma} {\rm \cite{EFR} (Graph Removal Lemma)}. Let\/ $H$ be a fixed graph on\/ $m$ vertices. If a graph\/ $G$ on\/ $n$ vertices contains\/ $o(n^m)$ copies of\/ $H$ $($as\/ $n\to\infty)$ then all copies can be destroyed by removing\/ $o(n^2)$ edges from\/ $G$.
\end{thm}

Let the pair $( E, \K )$ be a dominating set where $E\subset {[n]\choose 2}$ and $\K\subset {[n]\choose k}$.
We need to give a lower bound on $|E|+|\K|$. Since $( E, \K )$ is a dominating set, (i) and (ii) must hold.
For a given $e\in E$ there are exactly ${n-2\choose k-2}$  elements  $L\in {[n]\choose k}$ satisfying $E\subset L$.
Hence (i) implies
$${n-2\choose k-2}|E|+ |\K |\geq {n\choose k}.\eqno(6)$$
Similarly, counting the pairs in $[n]$ the inequality
$$|E|+ {k\choose 2}|\K |\geq {n\choose 2}\eqno(7)$$
is a consequence of (ii).
\begin{lem}\label{EKsimplelowerboundlem} If\/ $(E,\K)$ is a dominating set in\/ $G_{k,2}$ minimizing\/ $|E|+|\K|$, then we have
$$|E|\geq {1\over k(k-1)}n^2+o(n^2).$$
\end{lem}

{\noindent\bf Proof.} Because $E={[n]\choose 2}$ is a dominating set, an optimal dominating set satisfies
$|\K|=O(n^2)$. Dividing (6) by ${n-2\choose k-2}$ and using $k\geq 3$, we have
 $$|E|  \geq \frac{\binom{n}{k}}{\binom{n-2}{k-2}} - \frac{|\K|}{\binom{n-2}{k-2}} \geq \frac{n(n-1)}{k(k-1)} - \frac{O(n^2)(k-2)^{k-2}}{(ne)^{k-2}},$$
 which implies the statement of the lemma. \qed

\begin{lem}\label{EKlowerboundlem} If\/ $(E,\K)$ is a dominating set in\/ $G_{k,2}$ minimizing\/ $|E|+|\K|$, then we have
$$|E|\geq {1\over 2(k-1)}n^2+o(n^2).$$
\end{lem}

The reader might wonder why we formulated the weaker Lemma~\ref{EKsimplelowerboundlem}. The reason is that it will be used in the case
of $k=3$ in the proof of Lemma~\ref{EKlowerboundlem}.

\vskip 2mm

{\noindent\bf Proof.} Suppose that
$$|E|\leq {n\choose 2}{1\over k-1}-\left( {1\over k-1}-\varepsilon \right) n.$$
 Then the number of $k$-element subsets of $[n]$ containing no element of $E$ is at least
$$\left( {1\over k-1}-\varepsilon \right){n^{k-1}\over (k-1)^{k-2}}\eqno(8)$$
by Corollary~\ref{lovaszsimcor2}. Hence $|\K|$ must be at least this large. If this is more than the construction given in Section 2,
then a contradiction is obtained. However the inequality
$$\left( {1\over k-1}-\varepsilon \right){n^{k-1}\over (k-1)^{k-2}}>{k+3\over 2(k-1)(k+1)}n^2+o(n^2)\eqno(9)$$
is obvious for large $n$ if $k>3$. This contradiction shows the validity of the statement of the lemma when $k>3$.

When $k=3$, (9) is not true, so we do not arrive at a contradiction, but we know
$$|\K|\geq \left( {1\over 2}-\varepsilon \right){n^2\over 2}\eqno(10)$$
by (8).
Lemma~\ref{EKsimplelowerboundlem} and (10) lead to
$$|E|+|\K|\geq \left( {1\over 6}+\left( {1\over 4}-{\varepsilon \over 2} \right) \right)n^2+o(n^2).$$
This is more than ${3\over 8}n^2$ if $\varepsilon$ is small enough, contradicting the minimality of
the pair $(E,\K)$.  \qed

\vskip 2mm

It is easy to see that Lemma~\ref{EKlowerboundlem} and (7) give
$$|E|+|\K| \ge \frac{k^2+k-4}{k(k-1)^2}\,\frac{n^2}{2} + o(n^2),$$ but this is not strong enough.
We will improve (7). For $\h \subset {[n]\choose k}$, define the {\it 2-shadow} of $\h$ as
$\sigma_2 (\h)=\{ \{ i,j\}:\ i\not=j, \{i,j\}\subset H \mbox{ for some } H\in \h\}.$
Replacing $|\K|$ by $|\sigma_2 (\K)|$ is an essential improvement, since
every pair $e$ is counted only once even when it is contained in several members of $\K$. A further improvement
is obtained when the structure of $\K$ is also exploited.

For $\T\subset {[n]\choose k}$, denote by $Q(\T)$ denote the graph with vertex set $\T$, where $A,B\in \T$ $(A\not= B)$ are joined by an edge if $|A\cap B|\geq 2$.
Let $\K_0\subset \K$ be the set of those members which contain at least one element of $E$ as a subset.
Consider the components of $Q(\K -\K_0)$. Choose an integer $s>1$ and let $\K_2$ be the set of vertices
(of course, they are $k$-element subsets of $[n]$) belonging to a component of size at least $s$.
Finally, define $\K_1=\K -\K_2-\K_0$: this is the family of those members which contain no element of $E$ and are vertices of a component of $Q(\K-\K_0)$ of size at most $s-1$.  The improvement of (7) that will be really used is
$$|E|+\left( {k\choose 2}-1\right)|\K_0|+|\sigma_2 (\K_1)|+|\sigma_2 (\K_2)|\geq {n\choose 2}.\eqno(11)$$

The next two lemmas will show that the coefficient ${k\choose 2}$ in (7) can be replaced by a constant that is ``almost" ${k\choose 2}-1$ for members of $\K_2$.

\begin{lem}\label{firstshadowlem} If\/ $\T\subset {[n]\choose k}$,\/ $|\T|=t$ and\/ $Q(\T)$ is connected then\/ $\sigma_2(\T)\leq t\left( {k\choose 2}-1\right)+1$ holds.
\end{lem}

{\noindent\bf Proof.} We use induction on $t$. The case $t=1$ is trivial. Suppose that it is true for $t-1$ and prove for $t$.
A connected graph always has a vertex whose removal keeps the connectedness, for instance a leaf of a spanning tree.
Delete this vertex $T\in \T$. The graph $Q(\T- \{ T\})$ is connected and has $t-1$ vertices; by the induction hypothesis  $\sigma_2(\T-\{ T\})\leq (t-1)\left( {k\choose 2}-1\right)+1$. But $T$ gives at most
${k\choose 2}-1$ new elements to $\sigma (\T-\{ T\})$ since one of its pairs is a subset of a member of $\T-\{ T\}$, because of the connectivity of $Q(\T)$. \qed

\begin{lem}\label{secondshadowlem} The following inequality holds for $\K_2$:
$$|\sigma_2 (\K_2)|\leq \left( {k\choose 2}-1+{1\over s}\right) |\K_2|.$$
\end{lem}

{\noindent\bf Proof.} Let the components of $Q(\K_2)$ be $\T_1, \T_2, \ldots ,\T_r$ where $|\T_i|=t_i\geq s$.
Using Lemma~\ref{firstshadowlem}, we have
$$ |\sigma_2 (\K_2)|= \sum_{i=1}^r\sigma_2(\T_i)\leq \sum_{i=1}^rt_i \left( {k\choose 2}-1+{1\over t_i}\right)$$
$$\leq \sum_{i=1}^rt_i \left( {k\choose 2}-1+{1\over s}\right) =|\K_2|\left( {k\choose 2}-1+{1\over s}\right).$$
\qed

Now we show that $|\sigma_2 (\K_1)|$ is negligible.
\begin{lem}\label{shadowk1smalllem} The following inequality holds for $\K_1$:
$$|\sigma_2 (\K_1)|=o(n^2).$$
\end{lem}

{\noindent\bf Proof.} Let the vertex sets of the components of $Q(\K_1)$ be $\R_1, \R_2,\ldots ,\R_u$. Each vertex of each $\R_i$ determines a complete graph on $k$ vertices in $[n]$. Now let $\hat{Q} (\R_i)$ denote the graph obtained by taking the union of the edge sets of these complete graphs defined by the vertices of $\R_i$.

{\it Claim 1.} The total number $u$ of the components of $Q(\K_1)$ is $O(n^2)$.

{\it Proof.}  The graphs $\hat{Q} (\R_i)$ and $\hat{Q} (\R_j) (i\not=j)$ have no common edges therefore the total number of edges of the graphs $\hat{Q} (\R_i) (1\leq i\leq u)$ is at most $n\choose 2$. Hence the number of graphs
is also at most ${n \choose 2}$.

{\it Claim 2.} The number $f(n,k,s)$ of non-isomorphic copies $\hat{Q} (\R_i)$ is at most $2^{{ks\choose k}}$.

This bound is very weak, but we only need its independence of $n$. Here $\hat{Q} (\R_i)$ is a union of at most $s-1$
complete graphs on $k$ vertices. The total number of vertices is at most $sk$. (Actually it is only at most $k+(s-2)(k-2)$ but it is not important from our point of view.). On this underlying set of vertices one can choose a complete $k$-graph in at most ${sk\choose k}$ ways. To obtain $\hat{Q} (\R_i)$ we have to choose at most $s-1$ of the $k$-element sets. The number of choices is upper-bounded if any number of $k$-element sets is allowed to be chosen. The number of such choices is given in the Claim.

Return now to the proof of the lemma. Each $\hat{Q} (\R_i)$ is isomorphic to one of the graphs $H_1, H_2, \ldots , H_{f(n,k,s)}$. For a fixed $i$ the number of copies of $H_i$ is $O(n^2)$ by Claim 1. This is $o(n^{|H_i|})$ since $|H_i|\geq k\geq 3$. Therefore the Graph Removal Lemma can be applied: there are $o(n^2)$ edges in the graph $\cup\hat{Q}(\R_i)$  such that all copies of $H_i$ contain one of  them. Hence the number of copies of $H_i$ is also at most $o(n^2)$. The total number of components $\hat{Q} (\R_i)$  is at most $f(k,s)o(n^2)$. One component contains fewer than $s$ members of $\K_1$, hence $|\K_1|\leq sf(k,s)o(n^2)$, proving the lemma. \qed

\vskip 2mm

Now we can return to the proof of the lower bound in Theorem~\ref{gammagk2thm}. Inequal\-ity (11), Lemmas~\ref{secondshadowlem} and \ref{shadowk1smalllem} give
$$|E|+\left( {k\choose 2}-1+{1\over s}\right)|\K| \geq {n^2\over 2}+o(n^2).$$
Since $s$ can be arbitrarily large and $|\K|=O(n^2)$, this implies
$$|E|+\left( {k\choose 2}-1\right)|\K| \geq {n^2\over 2}+o(n^2).\eqno(12)$$
Multiplying (12) by ${1\over  {k\choose 2}- 1}$, the inequality in Lemma~\ref{EKlowerboundlem} by $1- {1\over  {k\choose 2}- 1}$ and adding them, the desired inequality is obtained. \qed

\section{Lower bound: short proof}
Similarly to the previous proof, let $(E,\K)$ be a dominating set in $G_{k,2}$ of minimum size. Consider the complement of $E$ on the vertex set $[n]$: $F={[n]\choose 2}-E$. A $k$-element subset $K\in {[n]\choose k}$ is {\it critical} if all of its ${k\choose 2}$ pairs are in $F$; denote the family of critical sets by $\K^*$. If $K\in \K^*$ then, by (i), $K\in \K$ must also hold. Hence we have $\K^*\subset \K$ and consequently $|\K^*|\leq |\K|$. Since $\binom{[n]}{2}$ is a dominating set of quadratic order and our dominating set is optimal, $|\K^*|\leq |\K|=O(n^2)$.
Now, by the Graph Removal Lemma we can find $o(n^2)$ edges in $F$ such that every member of $\K^*$ contains one of them
as a subset; denote by $H$ the set of these edges. Therefore we have
$$|H|=o(n^2).\eqno(13)$$
It is clear that
$$F-H {\mbox{ contains no complete graph on $k$ vertices}}.\eqno(14)$$

Consider now the slightly enlarged dominating set $(E\cup H, \K)$.
Here $E\cup H$ is the complement of $F-H$. Tur\'an's theorem can be applied for $F-H$ by (14). Therefore,
$|E\cup H|$ is at least as large as the number of edges of the graph obtained by $k-1$ vertex disjoint complete graphs
on $\lceil {n\over k-1} \rceil$ and $\lfloor {n\over k-1} \rfloor$ vertices, respectively. Hence we have 
$$|E|+ |H|\geq {n^2\over 2(k-1)}+o(n^2).\eqno(15)$$

By (14) every $K\in \K$ covers at most ${k\choose 2}-1$ ``new" edges, therefore (i) implies
$$|E|+|H|+ \left( {k\choose 2}-1\right) |\K|\geq {n\choose 2}.\eqno(16)$$
Taking (13) into account, (15) and (16) become
$$|E|\geq {n^2\over 2(k-1)}+o(n^2)\eqno(17)$$
and
$$|E|+\left( {k\choose 2}-1\right) |\K|\geq {n\choose 2} + o(n^2),\eqno(18)$$
respectively. Multiplying (18) by ${1\over  {k\choose 2}-1}$, (17) by $1- {1\over  {k\choose 2}-1}$, and adding them, the desired inequality is obtained.\qed

\section{The case $\ell >2$}

The case of $\ell >2$ is more complex.
A major difficulty is that the Tur\'an problem is still open for
 hypergraphs.
Nevertheless, should it become solved, the methods introduced above for
 $\gamma(G_{k,\ell})$ would almost surely be applicable for infinitely
 many pairs $k,\ell$.

More explicitly, the following can be done.
Denoting by $K_k^{(\ell)}$ the complete $\ell$-uniform hypergraph on $k$
 vertices ($k>\ell\geq 2$), let us write the Tur\'an numbers in the form
 $$
   \mathrm{ex}(n,K_k^{(\ell)}) = (1 - \alpha_{k,\ell})
     {n\choose \ell} + o(n^\ell) .
 $$
Hence, turning to complementary formulation,
 the minimum number of sets in an $\ell$-uniform set system $\F$ over
 $[n]$ such that every $k$-element set contains at least one $F\in\F$ is
 $(\alpha_{k,\ell}+o(1)) {n\choose \ell} + o(n^\ell)$ as $n\to\infty$.
By Tur\'an's theorem we have $\alpha_{k,2}= {1\over k-1}$, but
 the value of $\alpha_{k,\ell}$ is not known for $\ell\geq 3$.

Concerning the domination problem,
 as in the case of $\ell = 2$, we can start with a smallest
 dominating set $D$ of $G_{k,\ell}$.
Then, applying the Hypergraph Removal Lemma we can modify $D$ to a
 $D'$ which is still not much larger but already dominates
 the entire $[n]\choose k$ by the family $\cL\subset D'$ of
  selected $\ell$-element sets.
In this way the $k$-element sets in $D'$ are not needed for themselves
 in a set dominating the
  $k$-th level, hence the current condition
 on the sets to be selected from
 $[n]\choose k$ into $D'$ is that they should cover those members of
  $[n]\choose \ell$ which have not been selected into $D'$.
Since those $\ell$-element sets form a $K_k^{(\ell)}$-free family,
 each $k$-element subset can contain at most ${k\choose \ell}-1$ of them.
In this way we obtain the following general lower bound for every pair of fixed $ k > \ell > 2$.

 \begin{thm}\label{al-ellthm}
  For every fixed\/ $k>\ell\ge 2$, as\/ $n\rightarrow\infty$,
	$$\gamma(G_{k,\ell})\geq \left( \alpha_{k,\ell}+{1-\alpha_{k,
     \ell}\over {k\choose \ell}-1}\right) \!
     {n\choose \ell}  +o(n^{\ell}).$$
 \end{thm}

Some very interesting cases occur for particular values of $k$ and $\ell$,
 especially for $\ell=3$ and $k=4,5$.
We proceed in reverse order, since the situation with $k=4$ is more delicate.

\begin{prop} We have
 $\gamma(G_{5,3}) \leq \frac{1}{3}{n\choose 3} + o(n^3)$.
\end{prop}

{\noindent\bf Proof.}
A transparent construction showing the upper bound
 $\alpha_{5,3}\leq {1\over 4}$ is obtained by splitting the
 vertex set into two equal parts, say $A_1$ and $A_2$, and selecting
 all triplets inside each part.
These $\frac{1}{4}{n\choose 3} + o(n^3)$ 3-sets dominate all 5-sets.
It remains to find suitable 5-sets which cover all 3-sets
 meeting both $A_1$ and $A_2$.
The collection $\T$ of those 3-sets has size
 $\frac{3}{4}{n\choose 3} + o(n^3)$.

Consider the family $\F$ of all 5-sets $F$ with $|F\cap A_i|=2$ and
 $|F\cap A_{3-i}|=3$, for $i=1,2$.
Each $F\in \F$ contains 9 sets $T\in \T$.
Hence, this structure yields a 9-uniform hypergraph which is regular
 of degree $(c+o(1)) n^2$ for a $c>0$ if $n$ is even, and not far
  from being regular if $n$ is odd.
 On the other hand, the two-degrees are $O(n)$.
Thus the existence of an asymptotically optimal cover with
 $\frac{1}{12}{n\choose 3} + o(n^3)$ members of $\F$ follows by
 Theorem \ref{frthm}.
This yields $\gamma(G_{5,3}) \leq
 (\frac{1}{4} + \frac{1}{12}){n\choose 3} + o(n^3)$, as needed.
\qed

\begin{prop} We have
 $\gamma(G_{4,3}) \leq \frac{17}{27}{n\choose 3} + o(n^3)$.
\end{prop}

{\noindent\bf Proof.}
Assume for simplicity that $n$ is divisible by 3, and partition $[n]$
 into three sets $A_1,A_2,A_3$ of size $n/3$ each.
Number the elements from 1 to $n/3$ in $A_1$, from $n/3+1$ to $2n/3$ in $A_2$,
 and from $2n/3+1$ to $n$ in $A_3$.

We say that a set $F$ is of type $(p,q,r)$ if, for some $i\in\{1,2,3\}$,
 it satisfies $|F\cap A_{i}|=p$, $|F\cap A_{i+1}|=q$, $|F\cap A_{i+2}|=r$.
Subscript addition is taken modulo 3, i.e.\ the types
 $(p,q,r)$, $(q,r,p)$, $(r,p,q)$ mean exactly the same.

In the classical construction for the Tur\'an problem concerning
 $K_4^{(3)}$ one takes the 3-sets of types $(3,0,0)$ and $(2,1,0)$.
These dominate all 4-subsets of $[n]$.
Then we need to find 4-tuples of types $(2,1,1)$ and $(1,3,0)$ which
 dominate the non-selected 3-sets, namely those of
  types $(1,1,1)$ and $(1,2,0)$.
We now consider:
 \begin{itemize}
  \item all 4-sets of type $(1,3,0)$, and
  \item those 4-sets of type $(2,1,1)$ in which the sum of elements is even.
 \end{itemize}

A 3-set of type $(1,1,1)$ is not contained in any 4-set of type $(1,3,0)$,
 but it can be completed to a 4-set of type $(2,1,1)$ in $n/6-O(1)$
 different ways from each $A_i$, hence its degree is $n/2-O(1)$.

A 3-set of type $(1,2,0)$ can be completed to a 4-set of type $(1,3,0)$
 in $n/3-O(1)$ different ways from the ``middle'' part, and to a 4-set
 of type $(2,1,1)$ in $n/6-O(1)$ different ways from the ``last'' part,
 hence its degree is $n/3 + n/6 - O(n) = n/2-O(1)$.

It follows that a nearly regular 3-uniform hypergraph has been obtained
 whose vertices are the 3-sets of types $(1,1,1)$ and $(1,2,0)$,
 and its edges correspond to the 4-sets listed above.
Thus, Theorem \ref{frthm} implies that an asymptotically optimal cover exists.
Simple computation yields the claimed upper bound of
 $\frac{17}{27}{n\choose 3} + o(n^3)$.
\qed

\vskip 2mm

These considerations lead to the following very reasonable conjectures.

\begin{conj}\label{gammag53conj}
 $\gamma(G_{5,3}) = \frac{1}{3}{n\choose 3} + o(n^3)$ holds.
\end{conj}

\begin{conj}\label{gammag43conj}
 $\gamma(G_{4,3}) = \frac{17}{27}{n\choose 3} + o(n^3)$ holds.
\end{conj}

Many further interesting pairs $(k,\ell)$ can be considered;
 for instance, it is tempting to guess that the hypergraphs
 $K_n^{(3)}-3K_{n/3}^{(3)}$ are good candidates for providing the correct
 asymptotics of ex$(n, K_7^{(3)})$.
In cases where the analogous constructions are essentially tight, the above
 method is strong enough to determine the
 approximate value of $\gamma (G_{k,\ell})$.
In general, we formulate the following.

\begin{conj}\label{gammagklconj}
 The lower bound given in Theorem~\ref{al-ellthm} is asymptotically tight
  for all fixed\/ $k >\ell\ge 2$ as\/ $n\to\infty$.
\end{conj}

\section{Remarks}
\hskip 6.1mm {\bf 1.}$\enspace$Ervin Gy\H ori informed us that he also proved the case $k=3, \ell=2$ \cite{Gy}.

{\bf 2.}$\enspace$Y. Pandit, S.L. Sravanthi, S. Dara, and S.M. Hegde \cite{PSDH} considered the cases where $k>\lceil {n\over 2}\rceil$. Of course our previous conjecture and present theorem do not claim anything for this case.

{\bf 3.}$\enspace$Finally, let us call the attention to the somewhat related papers \cite{GHKLR} and \cite{GKLPPP}. In particular, by adapting Construction 13 of \cite{GKLPPP}, one can obtain a general upper bound for $\gamma(G_{k, \ell})$ which is a slight improvement over the trivial upper bound $\binom{n}{\ell}$.
\begin{prop}\label{simpleupperboundgammagklprop}
Let\/ $k$ and\/ $\ell$ be fixed, with\/ $k > \ell \ge 3$ and\/ $n\rightarrow \infty$. Then,
\[\gamma(G_{k, \ell}) \le \binom{n}{\ell} \!\! \left(1-\frac{k-\ell}{k-\ell+1}\left(1-\frac{1}{\ell}\right)^{\ell-1}+o(1)\right).\]
\end{prop}

\end{document}